\documentclass[reqno]{article}%
\usepackage{amsmath}
\usepackage{graphicx}
\usepackage{amsfonts}
\usepackage{amssymb}%
\newtheorem{theorem}{Theorem}

\newtheorem{conjecture}[theorem]{Conjecture}

\newtheorem{definition}[theorem]{Definition}

\newtheorem{lemma}[theorem]{Lemma}

\newenvironment{proof}[1][Proof]{\noindent{\textbf {#1}  }}  {\hfill$\Box$}

\begin{document}

\title{The Cycle-Complete graph Ramsey numbers}
\author{V. Nikiforov$^{\ref{f1}}$}
\maketitle

\begin{abstract}
In 1978 Erd\H{o}s, Faudree, Rousseau, and Schelp conjectured that
\[
r\left(  C_{p},K_{r}\right)  =\left(  p-1\right)  \left(  r-1\right)  +1.
\]
for every $p\geq r\geq3,$ except for $p=q=3.$ This has been proved for
$r\leq6,$ and for
\[
p\geq r^{2}-2r.
\]

In this note we prove the conjecture for $p\geq4r+2.$

\end{abstract}

\section{Introduction}

\footnotetext[1]{\label{f1} Department of Mathematical Sciences, University of
Memphis, Memphis, Tennessee, 38152.}

The problem of finding the Ramsey number $r\left(  C_{p},K_{r}\right)  $ has
attracted considerable attention in the last decades. In particular,
relatively good asymptotics are known in the case of $p$ fixed and $r$ large
(see \cite{Rou} for a recent survey). By contrast, the case $p\geq r$ is
poorly known. In \cite{BoEr} Bondy and Erd\H{o}s proved that for $r>3$ and
$p\geq r^{2}-2$ the following exact result holds%
\begin{equation}
r\left(  C_{p},K_{r}\right)  =\left(  p-1\right)  \left(  r-1\right)
+1.\label{maineq}%
\end{equation}
Later in \cite{EFRS} Erd\H{o}s, Faudree, Rousseau, and Schelp conjectured that
(\ref{maineq}) holds for every $p\geq r\geq3,$ except for $p=r=3.$ This has
been proved for $r=4$ in \cite{YHZ}, for $r=5$ in \cite{BJM}, and for $r=6$ in
\cite{Schr}. In \cite{Schr} Schiermeyer has also shown that (\ref{maineq})
holds for $r>3,$ $p\geq r^{2}-2r$.

In this note we prove that (\ref{maineq}) holds for all $r\geq3$ and
$p\geq4r+2.$

\section{Main results}

Our graph theoretic notation is standard (e.g., see \cite{Bol}). In
particular, for any vertex $u,$ $\Gamma\left(  u\right)  $ is the set of its
neighbors. A graph $G$ is said to be $H$\emph{-free} if it does not contain a
subgraph isomorphic to $H.$ A path with endvertices $u$ and $v$ is called an
$uv$-path; $P\left(  u,v\right)  $ denotes a path $P$ that joins $u$ to $v.$
We write $\left(  v_{1},...,v_{k},v_{1}\right)  $ for the cycle whose edges
are $\left(  v_{1},v_{2}\right)  ,...,\left(  v_{k-1},v_{k}\right)  ,\left(
v_{k},v_{1}\right)  .$ The set $\left\{  m,m+1,...,n\right\}  $ is denoted by
$\left[  m,n\right]  $ and $\left[  n\right]  $ is the set $\left\{
1,...,n\right\}  .$ An \emph{interval} of length $l$ is a set of $l\geq0$
consecutive integers. Given a graph $G$ and two distinct vertices $u,v$ of $G$
we denote by $R_{G}\left(  u,v\right)  $ the set of the orders of all
$uv$-paths; we shorten $R_{G}\left(  u,v\right)  $ to $R\left(  u,v\right)  $
when $G$ is implicit.

Our main goal in this note is to prove the following theorem.

\begin{theorem}
\label{mainth} If $r\geq4$ and $p\geq4r+2$ then $r\left(  C_{p},K_{r}\right)
=\left(  p-1\right)  \left(  r-1\right)  +1.$
\end{theorem}

To shorten the proof of the theorem we have distilled its main parts into
several lemmas that might be of independent interest. The proofs of the lemmas
are presented in section 2.5.

Our main tool will be a particular class of graphs containing a Hamiltonian
cycle together with a rich set of chords; we call such graphs \emph{saws} and
study them in subsection 2.3.

\subsection{General preliminary lemmas}

Following Burr and Erd\H{o}s \cite{BEFRS} we call a connected graph $H$
$r$\emph{-good} if the Ramsey number $r\left(  H,K_{r}\right)  $ of the pair
$\left(  H,K_{r}\right)  $ satisfies
\[
r\left(  H,K_{r}\right)  =\left(  r-1\right)  \left(  \left|  H\right|
-1\right)  +1.
\]

The following lemma could be regarded as a general result in the theory of the
$r$-good graphs.

\begin{lemma}
\label{le3} Suppose $H$ is a graph such that $r\left(  K_{s},H\right)  \leq
sp+1$ for every $s\leq r.$ Then every $H$-free graph $G$ of order $pr+1$ and
with $\alpha\left(  G\right)  \leq r$ is $2$-connected.
\end{lemma}

Erd\H{o}s and Gallai have proved in \cite{ErGa}, p. 345, the following assertion.

\begin{theorem}
[Erd\H{o}s and Gallai]\label{thErGa} If $G$ is a $2$-connected graph with
$d\left(  w\right)  \geq\delta$ for all $w\neq u,v,$ then there is a $uv$-path
of order at least $\delta+1$ in $G.$\hfill$\square$
\end{theorem}

This versatile result can be further extended in some particular cases; we
consider such extensions in the following two lemmas.

\begin{lemma}
\label{leEG} Let $G$ be a $2$-connected graph and $u,v$ be two vertices such
that $G-u-v$ is a union of two disjoint nonempty graphs $G_{1}$ and $G_{2}.$
If $d\left(  w\right)  \geq\delta,$ for every $w\neq u,v,$ then there is a
$uv$-path of order at least $\delta+1$ that has no vertices in common with
$G_{1}.$
\end{lemma}

\begin{lemma}
\label{ErGa1} Let $G$ be a $2$-connected graph and $x\in V\left(  G\right)  $.
If $d\left(  y\right)  \geq\delta,$ for every $y\neq x,$ then for every two
vertices $u$ and $v$ there is a $uv$-path of order at least $\delta+1$.
\end{lemma}

\subsection{The Chopping and the Collating Lemmas}

\begin{definition}
Let $P$ be a $uv$-path; a \textbf{reduction} of $P$ is a $uv$-path $Q$ such
that all vertices of $Q$ belong to $P.$ A $q$\textbf{-reduction} of $P$ is a
reduction of order $q.$
\end{definition}

The following two lemmas will be cornerstones in the proof of Theorem
\ref{mainth}; the first one is called the Chopping Lemma.

\begin{lemma}
\label{lechop} Let $G$ be a graph with $\alpha\left(  G\right)  \leq\alpha,$
and $x,y$ be two distinct vertices of $G,$ and $P$ be a $xy$-path of order
$l.$ Then, for every interval $I$ with
\[
I\subset\left[  l\right]  ,\text{ }\left\vert I\right\vert =2\alpha,
\]
there is a $q$-reduction of $P$ for some $q\in I$.
\end{lemma}

The following lemma summarizes the main tool in the proof of Theorem
\ref{mainth}; it is called the Collating Lemma.

\begin{lemma}
\label{lecoll} Suppose $G$ is a graph and $V\left(  G\right)  =V_{1}\cup
V_{2}$ is a nontrivial partition. Let $\left(  x_{1},x_{2}\right)  ,$ $\left(
y_{1},y_{2}\right)  $ be two disjoint edges such that $x_{1},y_{1}\in V_{1}$
and $x_{2},y_{2}\in V_{2}.$ Let $G_{1}=G\left[  V_{1}\right]  ,$
$G_{2}=G\left[  V_{2}\right]  ,$ and $a,b,k,l_{1},l_{2}$ be positive integers
such that:

\emph{(i) }$b-a\geq k-1;$

\emph{(ii)} $\left[  a,b\right]  \subset R_{G_{1}}\left(  x_{1},y_{1}\right)
;$

\emph{(iii)} for every interval $I\subset\left[  l_{1},l_{2}\right]  $ with
$\left|  I\right|  =k,$ $I\cap R_{G_{2}}\left(  x_{2},y_{2}\right)
\neq\varnothing.$

Then for every $s\in\left[  a+l_{1}+k,b+l_{2}\right]  $ there is a cycle of
order $s$ in $G.$
\end{lemma}

\subsection{Saws}

\begin{definition}
Let $S$ be a graph with $V\left(  S\right)  =\left\{  v_{1},...,v_{2k+1}%
\right\}  $; $S$ is called a \textbf{saw} if it contains the Hamiltonian cycle
$\left(  v_{1},...,v_{2k+1},v_{1}\right)  $ together with the chords $\left(
v_{2s-1},v_{2s+1}\right)  $ for every $s\in\left[  k\right]  $. The cycle
\[
\left(  v_{1},...,v_{2k+1},v_{1}\right)
\]
is called the \textbf{backbone} of $S,$ and the value
\[
d\left(  S\right)  =\min\left\{  d_{S}\left(  v_{2k}\right)  ,d_{S}\left(
v_{2k+1}\right)  \right\}
\]
is called the \textbf{degree} of $S$.
\end{definition}

We shall define saws by identifying their backbones. Although saws look quite
complicated, there are simple sufficient conditions for the existence of large
saws as shown in the following lemma.

\begin{lemma}
\label{le1} If $G$ is a graph with minimal degree $\delta\left(  G\right)
\geq p$ and independence number $\alpha\left(  G\right)  \leq r$ then $G$ has
a saw of degree at least $p-r$.
\end{lemma}

The following four Lemmas show that saws contain many paths of consecutive
lengths; actually, we introduce and study them exactly for that reason.

\begin{lemma}
\label{pr1} Let $S=\left(  v_{1},...,v_{2k+1},v_{1}\right)  $ be a saw, and
$P\left(  v_{i},v_{j}\right)  $ be a path along the cycle $\left(
v_{1},...,v_{2k+1},v_{1}\right)  $ of order $l.$ Then:

\emph{(i)} if the edge $\left(  v_{1},v_{2k+1}\right)  $ does not belong to
$P\left(  v_{i},v_{j}\right)  ,$ then $P\left(  v_{i},v_{j}\right)  $ has a
$q$-reduction for every%
\[
q\in\left[  \left\lceil \frac{l}{2}\right\rceil +1,l\right]  ;
\]

\emph{(ii)} if the edge $\left(  v_{1},v_{2k+1}\right)  $ belongs to $P\left(
v_{i},v_{j}\right)  ,$ then $P\left(  v_{i},v_{j}\right)  $ has a
$q$-reduction for every%
\[
q\in\left[  \left\lfloor \frac{l}{2}\right\rfloor +2,l\right]  ;
\]
In particular,\ if $P(x,y)$ is a path of order $2k+1$ along the cycle $\left(
v_{1},...v_{2k+1},v_{1}\right)  $ then $P$ has $q$-reductions for every
$q\in\left[  k+2,2k+1\right]  .$
\end{lemma}

Lemma \ref{pr1} implies that if $S=\left(  v_{1},...v_{2k+1},v_{1}\right)  $
is a saw then, for every two consecutive vertices $x,y$ along the cycle
$\left(  v_{1},...v_{2k+1},v_{1}\right)  ,$ we have
\[
\left[  k+2,2k+1\right]  \subset R_{S}\left(  x,y\right)  .
\]
It turns out that if the degree of $S$ is large compared to its order then,
for every two vertices $x,y,$ the set $R_{S}\left(  x,y\right)  $ contains
even larger intervals. We shall distinguish three different cases of pairs
$\left(  x,y\right)  \in V\left(  S\right)  ;$ each case is considered
separately in one of the following three lemmas.

\begin{lemma}
\label{super} Let $S=\left(  v_{1},...,v_{2k+1},v_{1}\right)  $ be a saw of
degree $d\left(  S\right)  \geq2\left(  2k+1\right)  /3.$ Then
\[
\left[  2,2k+1\right]  \subset R_{S}\left(  v_{2k},v_{2k+1}\right)  .
\]

\end{lemma}

\begin{lemma}
\label{lux} Let $S=\left(  v_{1},...,v_{2k+1},v_{1}\right)  $ be a saw with
$d\left(  S\right)  =d,$ and $x,y$ be two consecutive vertices along the cycle
$\left(  v_{1},...,v_{2k+1},v_{1}\right)  .$ Then
\[
\left[  2k-d+6,2k+1\right]  \subset R_{S}\left(  x,y\right)  .
\]

\end{lemma}

\begin{lemma}
\label{flat} Let $S=\left(  v_{1},...,v_{2k+1},v_{1}\right)  $ be a saw of
degree $d\left(  S\right)  =d\geq k,$ and $x,y$ be two distinct vertices of
$S.$ Then there exists some $l>d$ such that
\[
\left[  l-\left\lceil \frac{d}{2}\right\rceil +5,l\right]  \subset
R_{S}\left(  x,y\right)  .
\]

\end{lemma}

In the following lemma we combine path chopping with path reduction to prove
the existence of cycles of consecutive lengths in every saw with bounded
independence number.

\begin{lemma}
\label{le1.3} Let $k\geq3$ and $S$ be a saw of order $2k+1$ with independence
number $\alpha\left(  S\right)  \leq r.$ If $2r\leq k$ then $S$ contains a
cycle of order $q$ for every
\[
q\in\left[  4r,2k+1\right]  .
\]

\end{lemma}

\subsection{Proof of the theorem}

\begin{proof}
[Proof of Theorem \ref{mainth}]To slightly simplify the notation we shall
actually prove $r\left(  K_{r+1},C_{p+1}\right)  =pr+1,$ under the constraints
$r\geq2,$ $p\geq4r+5.$

Observe that the disjoint union of $r$ complete graphs of order $p$ is a
$C_{p+1}$-free graph of order $rp$ and has no independent set on $r+1$
vertices. Thus, for every $p$ and $r$ we have
\[
r\left(  K_{r+1},C_{p+1}\right)  \geq pr+1;
\]
so, all we have to prove is the inequality
\begin{equation}
r\left(  K_{r+1},C_{p+1}\right)  \leq pr+1. \label{mainin}%
\end{equation}

We shall use induction on $r;$ for $r\leq5$ (\ref{mainin}) follows from the
earlier results in \cite{YHZ}, \cite{BJM}, and \cite{Schr}; so we assume that
$r\geq6$ and (\ref{mainin}) holds for all $r^{\prime}<r.$ Assume
(\ref{mainin}) does not hold for $r$ and $p,$ and let $G$ be a $C_{p+1}$-free
graph of order $rp+1$ with $\alpha\left(  G\right)  \leq r;$ we shall show
that these assumptions lead to a contradiction.

First we shall prove that $\delta\left(  G\right)  \geq p.$ Indeed, let $u$ be
a vertex of minimal degree in $G$ and $V^{\prime}$ be the set of vertices that
are not adjacent to $u$ and are distinct from $u.$ Clearly,
\[
\alpha\left(  G\left[  V^{\prime}\right]  \right)  \leq r-1
\]
and $G\left[  V^{\prime}\right]  $ is $C_{p+1}$-free; thus, by the induction
hypothesis,
\[
\left|  V^{\prime}\right|  \leq p\left(  r-1\right)  ,
\]
and therefore,
\[
\delta\left(  G\right)  =pr-\left|  V^{\prime}\right|  \geq p.
\]

Applying Lemma \ref{le3} with $H=C_{p+1},$ we see that $G$ is $2$-connected.

Applying Lemma \ref{le1}, we can find in $G$ a saw $S=\left(  v_{1}%
,...v_{2k+1},v_{1}\right)  $ with
\[
d\left(  S\right)  =d\geq p-r.
\]
Since $G$ is $C_{p+1}$-free, Lemma \ref{le1.3} implies $2k+1\leq p,$ so we
have%
\begin{equation}
p-r\leq d<2k+1\leq p. \label{basin}%
\end{equation}

Set $G^{\ast}=G-S.$ Observe that from $\delta\left(  G\right)  \geq p$ it
follows that every vertex $w\in V\left(  S\right)  $ has at least $p-2k$
neighbors in $G^{\ast}.$

In the rest of the proof we construct a cycle of order $p+1$ under various
assumptions about the connectivity of $G^{\ast}$ and the edges of $E\left(
V\left(  S\right)  ,V\left(  G^{\ast}\right)  \right)  .$ To achieve this goal
we combine a $x_{1}y_{1}$-path in $S$ with a $x_{2}y_{2}$-path in $G^{\ast},$
where $\left(  x_{1},x_{2}\right)  ,$ $\left(  y_{1},y_{2}\right)  $ are
disjoint edges. Depending on the location of $x_{1},y_{1},$ we use Lemma
\ref{super}, \ref{lux}, or \ref{flat}, to find in $S$ sufficiently many
$x_{1}y_{1}$-paths of consecutive orders. On the other hand, we find in
$G^{\ast}$ a sufficiently long $x_{2}y_{2}$-path and use the Chopping Lemma
(Lemma \ref{lechop}) to show that $R_{G^{\ast}}\left(  x_{2},y_{2}\right)  $
hits any sufficiently short interval. Finally, we apply the Collating Lemma
(Lemma \ref{lecoll}) to show the existence of cycles of specified order,
including $C_{p+1},$ and thus, obtain a contradiction.

Let us now give the details of the proof. It is not hard to see that every
graph $G^{\ast}$ has one of the following properties:

- $G^{\ast}$ is $2$-connected;

- $G^{\ast}$ is not connected and all its components are $2$-connected;

- $G^{\ast}$ contains a connected component that is not $2$-connected (it may
be $G^{\ast}$ itself).

\textbf{Case 1. }$G^{\ast}$\emph{ is }$2$\emph{-connected.}

First we shall prove that $v_{2k}$ and $v_{2k+1}$ have two different
neighbors
\[
x_{2}\in V\left(  G^{\ast}\right)  \cap\Gamma\left(  v_{2k}\right)  ,\text{
}y_{2}\in V\left(  G^{\ast}\right)  \cap\Gamma\left(  v_{2k+1}\right)  .
\]

Indeed, this is clear if one of the vertices $v_{2k}$ and $v_{2k+1}$ has two
or more neighbors in $G^{\ast}.$ Otherwise, in view of (\ref{basin}), we see
that $p=2k+1,$ and $v_{2k},v_{2k+1}$ are joined to every vertex of $S.$ But
then, if $v_{2k}$ and $v_{2k+1}$ have a neighbor in common in $G^{\ast},$ we
immediately obtain a $C_{p+1},$ a contradiction. Thus, $v_{2k}$ and $v_{2k+1}$
have two different neighbors $x_{2},y_{2}\in V\left(  G^{\ast}\right)  $.

Next we shall show that there is a path $P_{2}\left(  x_{2},y_{2}\right)  $ of
order at least $p-2k$ in $G^{\ast}.$ We see this immediately if $p=2k+1,$
since $G^{\ast}$ is connected. If $p>2k+1,$ the assertion follows from
\[
\delta\left(  G^{\ast}\right)  \geq p-2k-1,
\]
and Lemma \ref{ErGa1} applied to $G^{\ast}$.

Let the order of $P_{2}\left(  x_{2},y_{2}\right)  $ be $l\geq p-2k.$ Hence,
from the Chopping Lemma, for every interval $I$ with%
\[
I\subset\left[  l\right]  ,\text{ }\left|  I\right|  =2r,
\]
there is a $q$-reduction of $P_{2}\left(  x_{2},y_{2}\right)  $ for some $q\in
I.$ On the other hand, from Lemma \ref{super},
\[
\left[  2,2k+1\right]  \subset R_{S}\left(  v_{2k},v_{2k+1}\right)  .
\]
Thus, as%
\[
2k\geq p-r\geq2r,
\]
from the Collating Lemma, we see that $G$ contains a cycle of order $s$ for
every
\[
s\in\left[  2r+2,l+2k+1\right]  ,
\]
and hence, from
\[
l+2k\geq p\geq2r+1,
\]
$G$ contains a $C_{p+1},$ a contradiction.

In the sequel we shall suppose that $G^{\ast}$ is decomposed into blocks and
edges. Clearly at least one of the endblocks of $G^{\ast}$ has independence
number at most $\left(  r+1\right)  /2;$ let $B$ be an endblock of $G^{\ast}$
with $\alpha\left(  B\right)  \leq\left(  r+1\right)  /2.$ The following case
will appear in several contexts, so we shall consider it separately.

\textbf{Case 2.}\emph{ There are two consecutive vertices along }$S$\emph{
that are joined to two distinct vertices of }$B.$

Let $x_{2},y_{2}\in V\left(  B\right)  $ and $x_{1},y_{1}$ be consecutive
vertices along $S$ such that
\[
x_{1}\in V\left(  S\right)  \cap\Gamma\left(  u_{1}\right)  ,\text{ }x_{2}\in
V\left(  S\right)  \cap\Gamma\left(  u_{2}\right)  ,\text{ }x_{1}\neq x_{2}.
\]
Observe that, except possibly for the cutvertex of $B$, for every vertex $w\in
V\left(  B\right)  $ we have
\[
d_{B}\left(  w\right)  \geq d_{G}\left(  w\right)  -\left|  S\right|  \geq
p-2k-1.
\]
From Lemma \ref{ErGa1}, every two vertices of $B$ are joined in $B$ by a path
$P_{2}\left(  x_{2},y_{2}\right)  $ of order%
\[
l\geq p-2k.
\]
Set $r_{1}=\alpha\left(  B\right)  ;$ from the choice of $B$ we have
\[
2r_{1}\leq r+1;
\]
the Chopping Lemma implies that, for every interval $I$ with
\[
I\subset\left[  l\right]  ,\text{ }\left|  I\right|  =2r_{1}%
\]
there is a $q$-reduction of $P_{2}\left(  x_{2},y_{2}\right)  $ for some $q\in
I.$ On the other hand, from Lemma \ref{lux},
\[
\left[  2k-d+6,2k+1\right]  \subset R_{S}\left(  x_{1},y_{1}\right)  .
\]
As, from (\ref{basin}), we have
\[
2k-d+6\leq p-d+5\leq r+5,
\]
we see that%
\[
\left[  r+5,2k+1\right]  \subset R_{S}\left(  x_{1},y_{1}\right)  ,
\]
and from (\ref{basin}),
\[
2k-r\geq p-2r\geq r.
\]
Applying the Collating Lemma, we find that $G$ contains a cycle of order $s$
for every
\[
s\in\left[  r+2r_{1}+5,l+2k+1\right]  .
\]
Hence, from
\[
l+2k\geq p>2r+6\geq r+2r_{1}+5,
\]
$G$ contains a $C_{p+1},$ a contradiction.

In the sequel we shall assume that for every endblock $B$ of $G^{\ast}$ with
\[
\alpha\left(  B\right)  \leq\left(  r+1\right)  /2,
\]
there are no consecutive vertices along $S$ that are joined to two distinct
vertices of $B.$

\textbf{Case 3. }$G^{\ast}$\emph{ is not connected and all its components are
}$2$\emph{-connected.}

Suppose $G^{\ast}$ is a union of disjoint connected components. Clearly, there
is a component $G_{1}$ of $G$ with $\alpha\left(  G_{1}\right)  \leq r/2.$
Select $u_{1}\in V\left(  G_{1}\right)  $ to have the maximum number of
neighbors in $S$ among the vertices of $G_{1};$ then for every $w\in V\left(
G_{1}\right)  ,$ $w\neq u_{1}$ we have
\[
d_{G_{1}}\left(  w\right)  \geq\frac{p}{2}.
\]
Indeed, assume there is some $w\in V\left(  G_{1}\right)  ,$ $w\neq u_{1}$
with $d_{G_{1}}\left(  w\right)  <p/2.$ Then%
\[
\left|  \Gamma\left(  u_{1}\right)  \cap V\left(  S\right)  \right|
\geq\left|  \Gamma\left(  w\right)  \cap V\left(  S\right)  \right|
=d_{G}\left(  w\right)  -d_{G_{1}}\left(  w\right)  >p-\frac{p}{2}=\frac{p}%
{2},
\]
and since $\left|  V\left(  S\right)  \right|  =2k+1\leq p,$ we see that
$u_{1}$ and $w$ are joined to two consecutive vertices along $S,$ a contradiction.

Since $G$ is $2$-connected, there are vertices $x_{1},y_{1}\in V\left(
S\right)  $ and $x_{2},y_{2}\in V\left(  B\right)  $ such that $\left(
x_{1},x_{2}\right)  $ and $\left(  y_{1},y_{2}\right)  $ are disjoint edges.
Applying Lemma \ref{ErGa1} to $G^{\ast}$, we see that there is a $x_{2}y_{2}%
$-path $P_{2}\left(  x_{2},y_{2}\right)  $ in $G_{1}$ of length $l>p/2$. Set
\[
V_{1}=V\left(  S\right)  ,\text{ }V_{2}=V\left(  G_{1}\right)  ,\text{ }%
\alpha\left(  G_{1}\right)  =r_{1};
\]
In view of $r_{1}\leq r/2,$\ the Chopping Lemma implies that for every
interval $I$ with%
\[
I\subset\left[  l\right]  ,\text{ }\left|  I\right|  =2r_{1},
\]
there is a $q$-reduction of $P_{2}\left(  x_{2},y_{2}\right)  $ for some $q\in
I.$ On the other hand, from Lemma \ref{flat}, there is some $l_{1}>d,$ such
that%
\[
\left[  l_{1}-\left\lfloor \frac{d}{2}\right\rfloor +5,l_{1}\right]  \subset
R_{S}\left(  x_{1},y_{1}\right)  .
\]
Applying the Collating Lemma we find that $G$ contains a cycle of order $s$
for every
\[
s\in\left[  l_{1}-\left\lfloor \frac{d}{2}\right\rfloor +5+2r_{1}%
,l_{1}+l\right]  .
\]
Hence, from
\[
l_{1}+l>\frac{p}{2}+d\geq\frac{p}{2}+p-r\geq p
\]
and
\[
p\geq2r+5\geq r+2r_{1}+5,
\]
$G$ contains a $C_{p+1},$ a contradiction.

\textbf{Case 4.}\emph{ }$G^{\ast}$\emph{ contains a connected component that
is not }$2$\emph{-connected.}

Clearly we can select an endblock $B$ of $G^{\ast}$ with
\[
\alpha\left(  B\right)  =r_{1}\leq\frac{r+1}{2};
\]
let $z$ be the only cutvertex of $B.$

Select $u_{1}\in V\left(  B-z\right)  $ to have the maximum number of
neighbors in $S$ among the vertices of $B-z,$ i.e.,%
\begin{equation}
\left|  \Gamma\left(  u_{1}\right)  \cap V\left(  S\right)  \right|
=\max_{w\in V\left(  B-z\right)  }\left|  \Gamma\left(  w\right)  \cap
V\left(  S\right)  \right|  . \label{maxd}%
\end{equation}

We shall show that for every $w\in V\left(  B-z\right)  ,$ $w\neq u_{1}$ we
have%
\begin{equation}
d_{B}\left(  w\right)  \geq\frac{p}{2}. \label{mind}%
\end{equation}

Indeed, assume there is some $w\in V\left(  B\right)  ,$ $w\neq u_{1}$ with
$d_{B}\left(  w\right)  <p/2.$ Then%
\[
\left|  \Gamma\left(  u_{1}\right)  \cap V\left(  S\right)  \right|
\geq\left|  \Gamma\left(  w\right)  \cap V\left(  S\right)  \right|
=d_{G}\left(  w\right)  -d_{B}\left(  w\right)  >p-\frac{p}{2}=\frac{p}{2},
\]
and since $\left|  V\left(  S\right)  \right|  =2k+1\leq p,$ we see that
$u_{1}$ and $w$ are joined to two consecutive vertices along $S,$ a contradiction.

Since $B$ is $2$-connected and (\ref{mind}) holds, from Theorem \ref{thErGa}
we see that the vertex $u_{1}$ is joined to $z$ by a path $P_{1}\left(
u_{1},z\right)  $ of order $l>p/2$.

Since $G^{\ast}\backslash B$ is nonempty, there is some $u_{2}\in V\left(
G^{\ast}\right)  \backslash V\left(  B\right)  $ that is joined to a vertex of
$S$ - otherwise $G-z$ is not connected, contradicting the fact that $G$ is $2$-connected.

\textbf{Case 4.1.}\emph{ There exist }$u_{2},x_{1},x_{2}$\emph{ with }%
\[
u_{2}\in V\left(  G^{\ast}\right)  \backslash V\left(  B\right)
,\emph{\ }x_{1}\in V\left(  S\right)  \cap\Gamma\left(  u_{1}\right)
,\emph{\ }x_{2}\in V\left(  S\right)  \cap\Gamma\left(  u_{2}\right)
,\emph{\ }x_{1}\neq x_{2}.
\]

Set $\alpha\left(  G^{\ast}\backslash B\right)  =r_{2};$ clearly, $r_{1}%
+r_{2}\leq r+1.$ Select the shortest $zu_{2}$-path $P_{2}\left(
z,u_{2}\right)  ;$ set
\[
l_{1}=\left|  P_{2}\left(  z,u_{2}\right)  \right|  .
\]
Since, except for the vertex $z,$ the path $P_{2}\left(  z,u_{2}\right)  $ is
entirely in $G^{\ast}\backslash B$ and
\[
\alpha\left(  G^{\ast}\backslash B\right)  =r_{2},
\]
the Chopping Lemma implies that
\[
2\leq l_{1}\leq2r_{2}+1.
\]
The concatenation
\[
Q=\left(  P_{1}\left(  u_{1},z\right)  ,P_{2}\left(  z,u_{2}\right)  \right)
\]
is a $u_{1}u_{2}$-path with
\[
\frac{p}{2}+2<\left|  Q\right|  =l+l_{1}-1\leq l+2r_{2}.
\]
Observe that the Chopping Lemma applied to the path $P_{1}\left(
u_{1},z\right)  ,$ implies that for every interval $I$ with
\[
I\subset\left[  l\right]  ,\text{ }\left|  I\right|  =2r_{1},
\]
there is a $q$-reduction of $P_{1}\left(  u_{1},z\right)  $ for some $q\in I.$
Therefore, for every interval $I$ with
\[
I\subset\left[  l_{1}+2r_{1},l+l_{1}-1\right]  ,\text{ }\left|  I\right|
=2r_{1},
\]
there is a $q$-reduction of $Q\left(  u_{1},u_{2}\right)  $ for some $q\in I.$

On the other hand, from Lemma \ref{flat}, there is some $l_{2}>d,$ such that%
\[
\left[  l_{2}-\left\lfloor \frac{d}{2}\right\rfloor +5,l_{2}\right]  \subset
R_{S}\left(  x_{1},x_{2}\right)  .
\]
Set
\[
V_{1}=V\left(  S\right)  ,\text{ }V_{2}=V\left(  G^{\ast}\right)  ;
\]
applying the Collating Lemma with the partition $V\left(  G\right)  =V_{1}\cup
V_{2},$ we find that $G$ contains a cycle of order $s$ for every%
\[
s\in\left[  l_{1}+l_{2}+2r_{1}-\left\lfloor \frac{d}{2}\right\rfloor
+5,l_{1}+l_{2}+l-1\right]  .
\]
Now from
\[
l_{1}\leq2r_{2}+1,\text{ }l_{2}\leq p-1,\text{ }\left\lfloor \frac{d}%
{2}\right\rfloor \geq\frac{p-r-1}{2},\text{ }r_{1}+r_{2}\leq r+1
\]
it follows that
\begin{align*}
l_{1}+l_{2}+2r_{1}-\left\lfloor \frac{d}{2}\right\rfloor +5  &  \leq
p+2r_{1}+2r_{2}-\frac{p-r-1}{2}+5\\
&  \leq p+r-\frac{p-r-1}{2}+6\leq p+1.
\end{align*}
On the other hand, from%
\[
l_{1}\geq2,\text{ }l_{2}>d\geq p-r,\text{ }l>\frac{p}{2}%
\]
it follows that%
\[
l_{1}+l_{2}+l-1>p-r+\frac{p}{2}+1>p.
\]
Hence, $G$ contains a $C_{p+1},$ a contradiction.

\textbf{Case 4.2.}\emph{ There are no }$u_{2},x_{1},x_{2}$\emph{ with }%
\[
u_{2}\in V\left(  G^{\ast}\right)  \backslash V\left(  B\right)
,\emph{\ }x_{1}\in V\left(  S\right)  \cap\Gamma\left(  u_{1}\right)
,\emph{\ }x_{2}\in V\left(  S\right)  \cap\Gamma\left(  u_{2}\right)
,\emph{\ }x_{1}\neq x_{2}.
\]

Clearly the assumption implies that $u_{1}$ has exactly one neighbor
\[
x_{1}\in V\left(  S\right)  \cap\Gamma\left(  u_{1}\right)  ,
\]
and for every $w\in V\left(  G^{\ast}\right)  \backslash V\left(  B\right)  ,$
either $w$ has no neighbors in $V\left(  S\right)  ,$ or is joined exactly to
$x_{1}.$ Observe that the choice of $u_{1}$ implies that every $w\in V\left(
B-z\right)  $ has at most one neighbor in $V\left(  S\right)  .$ Therefore,
for every $w\in V\left(  B-z\right)  ,$
\[
d_{B}\left(  w\right)  \geq p-1,
\]
and from Lemma \ref{ErGa1}, every two vertices of $B$ are joined by a path of
order at least $p.$

On the other hand, there must be some $u_{2}\in V\left(  B\right)  ,$ distinct
from $u_{1}$ and having a neighbor%
\[
x_{2}\in V\left(  S\right)  \cap\Gamma\left(  u_{2}\right)  \text{, }x_{2}\neq
x_{1},
\]
otherwise $G-x_{1}$ is not connected, contradicting the fact that $G$ is
$2$-connected. Select
\[
u_{2}\in V\left(  B\right)  ,\text{ }x_{2}\in V\left(  S\right)  \cap
\Gamma\left(  u_{2}\right)  ,\text{ }x_{2}\neq x_{1}.
\]
We know that there is a path $P_{1}\left(  u_{1},u_{2}\right)  $ in $B$ of
order $l\geq p;$ The Chopping Lemma implies that for every interval $I$ with
\[
I\subset\left[  l\right]  ,\text{ }\left|  I\right|  =2r_{1},
\]
there is a $q$-reduction of $P_{1}\left(  u_{1},z\right)  $ for some $q\in I.$

On the other hand, from Lemma \ref{flat}, there is some $l_{2}>d,$ such that%
\[
\left[  l_{2}-\left\lfloor \frac{d}{2}\right\rfloor +5,l_{2}\right]  \subset
R_{S}\left(  x_{1},x_{2}\right)  .
\]
Exactly as in the previous case, we see that $G$ contains a $C_{p+1},$ a
contradiction, completing the proof.
\end{proof}

\subsection{Proofs of the lemmas}

\begin{proof}
[\textbf{Proof of Lemma \ref{le3}}]Assume first that $G$ is disconnected; say
$G$ is the union of two disjoint nonempty graphs $G_{1}$ and $G_{2}$. Since
both $G_{1}$ and $G_{2}$ are $H$-free and
\[
\alpha\left(  G\right)  =\alpha\left(  G_{1}\right)  +\alpha\left(
G_{2}\right)  ,
\]
we have%
\[
v\left(  G\right)  =v\left(  G_{1}\right)  +v\left(  G_{2}\right)  \leq
\alpha\left(  G_{1}\right)  p+\alpha\left(  G_{2}\right)  p\leq pr,
\]
a contradiction. Thus, $G$ is connected.

Assume now that $G$ is not $2$-connected and let $u$ be a cutvertex of $G.$
Then $G-u$ is the union of two disjoint graphs $G_{1}$ and $G_{2}$. Since both
$G_{1}$ and $G_{2}$ are $H$-free, and
\[
\alpha\left(  G_{1}\right)  +\alpha\left(  G_{2}\right)  \leq\alpha\left(
G\right)  ,
\]
we have%
\[
pr+1=v\left(  G\right)  =v\left(  G_{1}\right)  +v\left(  G_{2}\right)
+1\leq\alpha\left(  G_{1}\right)  p+\alpha\left(  G_{2}\right)  p+1\leq pr+1.
\]
Hence,
\[
v\left(  G_{1}\right)  =\alpha\left(  G_{1}\right)  p,\text{ }v\left(
G_{2}\right)  =\alpha\left(  G_{2}\right)  p
\]
and
\[
\alpha\left(  G\right)  =\alpha\left(  G_{1}\right)  +\alpha\left(
G_{2}\right)  .
\]
By the condition of the Lemma,
\[
\alpha\left(  G_{1}+u\right)  \geq\alpha\left(  G_{1}\right)  +1,\text{
}\alpha\left(  G_{2}+u\right)  \geq\alpha\left(  G_{2}\right)  +1;
\]
thus,
\[
\alpha\left(  G\right)  \geq\alpha\left(  G_{1}+u\right)  +\alpha\left(
G_{2}+u\right)  -1\geq r+1,
\]
a contradiction. Hence, $G$ is $2$-connected.
\end{proof}

\begin{proof}
[\textbf{Proof of Lemma \ref{leEG}}]Remove $G_{1},$ take a second copy of the
remaining graph and identify the vertices $u$ and $v$ in both copies. The
resulting graph satisfies the hypothesis of Theorem \ref{thErGa}, and
therefore, contains a $uv$-path of length at least $\delta+1.$ This path is
contained in entirely in $G$ and consequently has no vertices in common with
$G_{1}$.
\end{proof}

\begin{proof}
[\textbf{Proof of Lemma \ref{ErGa1}}]If $u=x$ or $v=x$ the assertion follows
from Theorem \ref{thErGa}, so assume $u\neq x$ and $v\neq x.$ Set $G^{\ast
}=G-x;$ clearly we have $\delta\left(  G^{\ast}\right)  \geq\delta-1.$ From
the $2$-connectivity of $G$ it follows that $G^{\ast}$ is connected.

\textbf{Case 1.}\emph{ }$G^{\ast}$\emph{ is not }$2$\emph{-connected.}

Let $y$ be a vertex such that $G^{\ast}-y$ is a union of two disjoint nonempty
graphs $G_{1}$ and $G_{2}.$ Then, from Lemma \ref{leEG}, there exist two paths
$P\left(  x,y\right)  $ and $Q\left(  x,y\right)  $ such that
\[
P\cap V\left(  G_{1}\right)  =\varnothing,\text{ }Q\cap V\left(  G_{2}\right)
=\varnothing,\text{ }\left|  P\right|  \geq\delta+1,\text{ }\left|  Q\right|
\geq\delta+1
\]
Clearly, the concatenation of $P$ and $Q$ is a cycle $C$ of order at least
$2\delta.$ Since $G$ is $2$-connected, there are two vertices $u_{1}$ and
$v_{1}$ of $C$ and two disjoint paths $P_{u}\left(  u,u_{1}\right)  $ and
$P_{v}\left(  v,v_{1}\right)  ,$ possibly of order 1. Select $Q\left(
u_{1},v_{1}\right)  $ to be the path that $u_{1}$ and $v_{1}$ cut from $C$
with $\left|  Q\right|  \geq\left|  C\right|  /2;$ the concatenation $\left(
P_{u},Q,P_{v}\right)  $ is a $uv$-path of order at least $\delta+1.$

\textbf{Case 2.}\emph{ }$G^{\ast}$\emph{ is }$2$\emph{-connected. }

From a theorem of Dirac , since $G^{\ast}$ is $2$-connected and $\delta\left(
H\right)  \geq\delta-1,$ there is a cycle of order at least $2\delta-2$ unless
$v\left(  G^{\ast}\right)  <2\delta-2$. Consider first the latter case.

\textbf{Case 2.1.}\emph{ }$v\left(  G^{\ast}\right)  <2\delta-2.$

As Bondy proved (attributing the result to Erd\H{o}s and Gallai) in
\cite{Bon}, Corollary 2.13, the assumption $v\left(  H\right)  <2\delta-2$,
together with $\delta\left(  G^{\ast}\right)  \geq\delta-1,$ implies that
$G^{\ast}$ is Hamilton-connected, i.e., every two vertices are connected by a
Hamiltonian path. Hence, we obtain a $uv$-path of order at least $\delta+1,$
unless $G^{\ast}$ is a complete graph of order $\delta.$ In the latter case
$x$ must be adjacent to every vertex of $G^{\ast},$ as, otherwise, there is
some $w\in V\left(  G^{\ast}\right)  $ with $d_{G}\left(  w\right)
=\delta-1,$ a contradiction.\ Now we trivially obtain a $uv$-path of order
$\delta+1.$

\textbf{Case 2.2.}\emph{ }$v\left(  G^{\ast}\right)  \geq2\delta-2.$

Hence, $G^{\ast}$ has a cycle $C$ of order at least $2\delta-2.$ From the
$2$-connectivity of $G^{\ast}$ it follows that there are two disjoint paths
$P_{1}\left(  u,u_{1}\right)  $ and $P_{2}\left(  v_{1},v\right)  ,$ where
$u_{1},v_{1}\in C.$ Select $Q\left(  u_{1},v_{1}\right)  $ to be the path that
$u_{1}$ and $v_{1}$ cut from $C$ with $\left|  Q\right|  \geq\left|  C\right|
/2;$ the concatenation
\[
\left(  P_{1}\left(  u,u_{1}\right)  ,Q\left(  u_{1},v_{1}\right)
,P_{2}\left(  v_{1},v\right)  \right)
\]
is a $uv$-path of order at least $\delta+1,$ unless the order of $C$ is
precisely $2\delta-2,$ $u,v\in C$ and the distance along $C$ between $u$ and
$v$ is $\delta-1.$ Let
\[
P_{1}=\left(  u,u_{1},...,u_{\delta-2},v\right)  ,\text{ }P_{2}=\left(
v,v_{1},...,v_{\delta-2},u\right)
\]
be the two paths joining $u$ and $v$ along $C,$ i.e. $C=\left(  P_{1}%
,P_{2}\right)  .$ Consider first the case of connected $G-u-v.$

\textbf{Case 2.2.1.}\emph{ }$G-u-v$\emph{ is connected.}

Hence, there is a path $P\left(  u_{i},v_{j}\right)  $ joining some $u_{i}%
\in\left\{  u_{1},...,u_{\delta-2}\right\}  $ to some $v_{j}\in\left\{
v_{1},...,v_{\delta-2}\right\}  $ and such that $P\left(  u_{1},v_{1}\right)
$ has no internal vertices in common with $C.$ Let $Q\left(  v_{j}%
,u_{i}\right)  $ be the path $P\left(  u_{i},v_{j}\right)  $ taken in reverse
order; we easily see that at least one of the paths
\begin{align*}
&  \left(  u,u_{1},...,u_{i-1},P\left(  u_{i},v_{j}\right)  ,v_{j+1}%
,,...,v_{\delta-2},v\right) \\
&  \left(  v,v_{1},...,v_{j-1},Q\left(  v_{j},u_{i}\right)  ,u_{i+1}%
,,...,u_{\delta-2},u\right)
\end{align*}

has order at least $\delta+1.$

\textbf{Case 2.2.2.}\emph{ }$G-u-v$\emph{ is disconnected.}

Let $G-u-v$ be the union of two disjoint nonempty graphs $G_{1}$ and $G_{2}.$
Without loss of generality we may suppose $x\in G_{2}.$ Since the graph
$G^{\ast}=G-G_{2}$ is $2$-connected and $d_{G^{\ast}}\left(  w\right)
\geq\delta$ for all $w\neq u,v$ then, from Theorem \ref{thErGa}, there is a
$uv$-path of order at least $\delta+1$ in $G^{\ast}$ and the proof is completed.
\end{proof}

\begin{proof}
[\textbf{Proof of Lemma \ref{lechop}}]If $l\leq2\alpha$ the assertion is
trivially true, so suppose $l\geq2\alpha+1.$ Observe that $G$ has no induced
path on $2\alpha+1$ vertices - otherwise choosing every other vertex along
such path we obtain an independent set on $r+1$ vertices. Hence, the first
$2\alpha+1$ vertices of $P$ induce a chord and there is a $q$-reduction
$P_{1}$ of $P$ for some
\[
q\in\left[  l-2\alpha,l-1\right]  .
\]
Setting $P_{0}=P$ and repeating the same argument as long as possible, we
obtain a sequence
\[
P_{0},P_{1},...,P_{s}%
\]
of reductions of $P$ such that for every $i=0,...,s-1,$%
\[
\left|  P_{i}\right|  \geq2\alpha+1,\text{ }\left|  P_{i}\right|  -2\alpha
\leq\left|  P_{i+1}\right|  <\left|  P_{i}\right|
\]
and $\left|  P_{s}\right|  \leq2\alpha.$ Clearly, every interval $I$ of length
$2\alpha$ in $\left[  l\right]  $ contains the order of some $\left|
P_{i}\right|  $ and the proof is completed.
\end{proof}

\begin{proof}
[\textbf{Proof of Lemma \ref{lecoll}}]Let $\left\{  s_{1},s_{2},...,s_{t}%
\right\}  =R_{G_{2}}\left(  x_{2},y_{2}\right)  \cap\left[  l_{1}%
,l_{2}\right]  $ and suppose
\[
l_{1}=s_{1}<s_{2}<...<s_{t}=l_{2}.
\]
From \emph{(iii)} we see that $s_{i+1}-s_{i}\leq k\leq b-a+1.$ Combining a
fixed path $P\left(  x_{2},y_{2}\right)  $ of order $s_{i}$ in $G_{2}$ with a
path $Q\left(  y_{1},x_{1}\right)  $ in $G_{1}$ of order $q$ for every
$q\in\left[  a,b\right]  $ we obtain a cycle of order $s$ for every
$s\in\left[  s_{i}+a,s_{i}+b\right]  .$ Since for every $i\in\left[
t-1\right]  $ we have
\[
s_{i+1}-s_{i}\leq b-a+1,
\]
the intervals
\[
\left[  s_{i}+a,s_{i}+b\right]  \text{, }\left[  s_{i+1}+a,s_{i+1}+b\right]
\]
are contiguous or overlap and the assertion follows.
\end{proof}

\begin{proof}
[\textbf{Proof of Lemma \ref{le1}}]Let $P=\left(  v_{1},...,v_{2t+1}\right)  $
be the longest path in $G$ such that the chord $\left(  v_{2s-1}%
,v_{2s+1}\right)  $ exists for every $s=1,...,t.$ The vertex $v_{2t+1}$ is
joined to at least $p-r$ of the vertices $v_{1},...,v_{2t}$ - otherwise the
set
\[
N=\Gamma\left(  v_{2t+1}\right)  \backslash\left\{  v_{1},...,v_{2t}\right\}
\]
contains at least
\[
d\left(  v_{2t+1}\right)  -\left(  p-r\right)  +1\geq r+1
\]
vertices and thus, $N$ induces an edge that extends $P$ by two more vertices.
By symmetry, $v_{2t}$ is joined to at least $p-r$ of the vertices
$v_{1},...v_{2t-1},v_{2t+1}.$ Let $i$ be the minimal index such that $v_{i}\in
P$ is joined to either $v_{2t}$ or $v_{2t+1};$ without loss of generality we
may assume that $v_{i}$ is joined to $v_{2t+1}.$ If $i$ is odd then the graph
induced by
\[
\left\{  v_{i},v_{i+1},...,v_{2t+1}\right\}
\]
is a saw of degree at least $p-r.$ If $i$ is even then the graph induced by
\[
\left\{  v_{i},v_{i-1},v_{i+1},...,v_{2t+1}\right\}
\]
is a saw of degree at least $p-r.$
\end{proof}

\begin{proof}
[\textbf{Proof of Lemma \ref{pr1}}]Observe that every $3$-path $\left(
v_{2s-1},v_{2s},v_{2s+1}\right)  $ along $P\left(  v_{i},v_{j}\right)  $ can
be replaced by the chord $\left(  v_{2s-1},v_{2s+1}\right)  $ shortening $P$
by $1$. Such a replacement can be done as many times as there are $3$-paths
$\left(  v_{2s-1},v_{2s},v_{2s+1}\right)  $ along $P\left(  v_{i}%
,v_{j}\right)  ,$ so, all we have to do is to estimate their number.

\textbf{Case (i). }\emph{ The edge }$\left(  v_{1},v_{2k+1}\right)  $\emph{
does not belong to }$P\left(  v_{i},v_{j}\right)  .$

In this case we have $j-i=l-1$ and
\[
P\left(  v_{i},v_{j}\right)  =\left(  v_{i},,v_{i+1},...,v_{i+l-1}\right)  .
\]
The number of the $3$-paths $\left(  v_{2s-1},v_{2s},v_{2s+1}\right)  $ along
$P\left(  v_{i},v_{j}\right)  $ is exactly the number of all $s$ such that
\[
i\leq2s-1<i+l-1,
\]
and it is at least $\left\lfloor l/2\right\rfloor -1$. Hence the assertion follows.

\textbf{Case (ii).}\emph{ The edge }$\left(  v_{1},v_{2k+1}\right)  $\emph{
belongs to }$P\left(  v_{i},v_{j}\right)  .$

In this case we have $j-i=l-2k-2,$ and
\[
P\left(  v_{i},v_{j}\right)  =\left(  v_{i},...,v_{2k+1},v_{1}%
,...,v_{i+l-2k-2}\right)  .
\]
The number of the $3$-paths $\left(  v_{2s-1},v_{2s},v_{2s+1}\right)  $ along
the path $\left(  v_{i},...,v_{2k+1}\right)  $ is
\[
\left\lfloor \frac{2k+1-i}{2}\right\rfloor ,
\]
and the number of the $3$-paths $\left(  v_{2s-1},v_{2s},v_{2s+1}\right)  $
along $\left(  v_{1},...,v_{i+l-2k-2}\right)  $ is%
\[
\left\lfloor \frac{i+l-2k-3}{2}\right\rfloor .
\]
Thus, the number of the $3$-paths $\left(  v_{2s-1},v_{2s},v_{2s+1}\right)  $
along $P\left(  v_{i},v_{j}\right)  $ is
\[
\left\lfloor \frac{2k+1-i}{2}\right\rfloor +\left\lfloor \frac{i+l-2k-3}%
{2}\right\rfloor \geq\left\lfloor \frac{l}{2}\right\rfloor -2
\]
and the assertion follows.
\end{proof}

\begin{proof}
[\textbf{Proof of Lemma \ref{super}}]From Lemma \ref{pr1} we have
\[
\left[  k+2,2k+1\right]  \subset R_{S}\left(  v_{2k},v_{2k+1}\right)  ,
\]
so we need to prove only
\[
\left[  2,k+2\right]  \subset R_{S}\left(  v_{2k},v_{2k+1}\right)  .
\]
In fact we shall prove the following more general assertion, implying the
required result:

\emph{Let }$G$\emph{ be a Hamiltonian graph of order }$n\geq5$\emph{ and let
}$\left(  v_{1},...,v_{n},v_{1}\right)  $\emph{ be a Hamiltonian cycle in
}$G.$\emph{ If }%
\[
d\left(  v_{1}\right)  +d\left(  v_{n}\right)  \geq\left(  4n-1\right)  /3
\]
\emph{then}%
\[
\left[  2,\left\lceil n/2\right\rceil +2\right]  \subset R\left(  v_{1}%
,v_{n}\right)  .
\]
Indeed, choose some $q\in\left[  2,\left\lceil n/2\right\rceil \right]  .$ Our
first goal is to find two vertices
\[
v_{i},v_{i+q}\in\left\{  v_{2},...,v_{n-1}\right\}
\]
such that
\[
e\left(  \left\{  v_{i},v_{i+q}\right\}  ,\left\{  v_{1},v_{n}\right\}
\right)  \geq3.
\]
Assume this assertion is not true and consider first the case $q>\left(
n-2\right)  /3.$ The pairs
\[
\left(  v_{2},v_{q+2}\right)  ,\left(  v_{3},v_{q+3}\right)  ,...,\left(
v_{q+1},v_{2q+1}\right)
\]
are disjoint and their union is the set $\left\{  v_{2},...,v_{2q+1}\right\}
.$ Hence, we have
\[
e\left(  \left\{  v_{1},v_{n}\right\}  ,\left\{  v_{2},...,v_{2q+1}\right\}
\right)  \leq2q,
\]
and thus,%
\begin{align*}
d\left(  v_{1}\right)  +d\left(  v_{n}\right)   &  =e\left(  \left\{
v_{1},v_{n}\right\}  ,\left\{  v_{2},...,v_{2q+1}\right\}  \right)  +e\left(
\left\{  v_{1},v_{n}\right\}  ,\left\{  v_{2q+2},...,v_{n-1}\right\}  \right)
+2\\
&  \leq2q+2\left(  n-2-2q\right)  +2=2n-2q-2\\
&  <2n-2-\frac{2\left(  n-2\right)  }{3}=\frac{4n-2}{3},
\end{align*}

a contradiction.

Let now $q\leq\left(  n-2\right)  /3$ and suppose
\[
n-2=qs+t,\left(  0\leq t\leq q-1\right)  .
\]
It is not hard to find a set of $\left\lceil qs/2\right\rceil $ disjoint pairs
of vertices $\left\{  v_{i},v_{i+q}\right\}  $ in $\left\{  v_{2}%
,...,v_{n-1}\right\}  .$ Since for every pair $\left\{  v_{i},v_{i+q}\right\}
$ we have by assumption
\[
e\left(  \left\{  v_{i},v_{i+q}\right\}  ,\left\{  v_{1},v_{n}\right\}
\right)  \leq2,
\]
we find that
\begin{align*}
d\left(  v_{1}\right)  +d\left(  v_{n}\right)   &  \leq2\left\lceil \frac
{qs}{2}\right\rceil +2\left(  n-2-2\left\lceil \frac{qs}{2}\right\rceil
\right)  +2=2n-2-2\left\lceil \frac{qs}{2}\right\rceil \\
&  \leq2n-2-\left(  qs-1\right)  =n+\left(  t+1\right)  \leq n+q\leq
\frac{4n-2}{3},
\end{align*}

a contradiction.

Therefore, there are two vertices $v_{i},v_{i+q}\in\left\{  v_{2}%
,...,v_{n-1}\right\}  $ such that
\[
e\left(  \left\{  v_{i},v_{i+q}\right\}  ,\left\{  v_{1},v_{n}\right\}
\right)  \geq3.
\]
Hence, either the edges $\left(  v_{1},v_{i}\right)  ,\left(  v_{n}%
,v_{i+q}\right)  ,$ or the edges $\left(  v_{1},v_{i+q}\right)  ,\left(
v_{n},v_{i}\right)  $ exist. So, either the path
\[
\left(  v_{1},v_{i},v_{i+1},...,v_{i+q},v_{n}\right)
\]
or the path
\[
\left(  v_{1},v_{i+q},v_{i+q-1},...,v_{i},v_{n}\right)
\]
exists, and we see that $q+2\in R\left(  v_{1},v_{n}\right)  .$ Hence,
\[
\left[  4,\left\lceil n/2\right\rceil +2\right]  \subset R\left(  v_{1}%
,v_{n}\right)  ,
\]
and since, obviously
\[
2\in R\left(  v_{1},v_{n}\right)  ,\text{ }3\in R\left(  v_{1},v_{n}\right)
,
\]
the proof is completed.
\end{proof}

\begin{proof}
[\textbf{Proof of Lemma \ref{lux}}]Set
\begin{equation}
t=2k-d \label{eq1}%
\end{equation}
and observe that the set
\[
M=\left\{  v_{1},...v_{2k}\right\}  \backslash\Gamma\left(  v_{2k+1}\right)
\]
has at most $t$ members. We assume that $\left\{  x,y\right\}  \neq\left\{
v_{2k},v_{2k+1}\right\}  $ since the case $\left\{  x,y\right\}  =\left\{
v_{2k},v_{2k+1}\right\}  $ is covered by Lemma \ref{super}. Thus, up to
labeling, there are only two different cases
\[
\left\{  x,y\right\}  \subset\left\{  v_{1},...,v_{2k}\right\}  ,
\]
and
\[
x=v_{1},\emph{\ }y=v_{2k+1}.
\]

\textbf{Case 1.}\emph{ }$\left\{  x,y\right\}  \subset\left\{  v_{1}%
,...,v_{2k}\right\}  $

Let $x=v_{j},$ $y=v_{j+1}$ and $PR\left(  j\right)  $ be the set of all pairs
of vertices $\left(  v_{i},v_{l}\right)  $ such that%
\begin{equation}
v_{i},v_{l}\in\Gamma\left(  v_{2k+1}\right)  ,\ 1\leq i\leq j<l\leq2k.
\label{pairs}%
\end{equation}
For every $\left(  v_{i},v_{l}\right)  \in PR\left(  j\right)  $ the value
$\left(  l-i\right)  $ is called its \emph{span}. Observe that if $\left(
v_{i},v_{l}\right)  \in PR\left(  j\right)  $ then the sequence%
\[
\left(  v_{j+1},...,v_{l},v_{2k+1},v_{i},...,v_{j}\right)
\]
is a $v_{j}v_{j+1}$-path of order $\left(  l-i+2\right)  $ and this motivates
the investigation of $PR\left(  j\right)  $ that follows.

Suppose $\left(  v_{h},v_{m}\right)  ,\left(  v_{i},v_{l}\right)  \in
PR\left(  j\right)  $ are distinct; we write
\[
\left(  v_{h},v_{m}\right)  \succ\left(  v_{i},v_{l}\right)
\]
if
\[
i\geq h\text{, and }l\leq m.
\]

We shall construct a sequence of $\left(  v_{i_{h}},v_{l_{h}}\right)  \in
PR\left(  j\right)  $ in the following way. Note first that $v_{2k+1}$ is
joined to both $v_{1}$ and $v_{2k}$ and thus $PR\left(  j\right)
\neq\varnothing$. Set $i_{1}=1,$ $l_{1}=2k.$

It turns out that if $\left(  v_{i_{h}},v_{l_{h}}\right)  \in PR\left(
j\right)  $ has a large span then there exists $\left(  v_{i_{h+1}}%
,v_{l_{h+1}}\right)  \in PR\left(  j\right)  $ such that
\[
\left(  v_{i_{h}},v_{l_{h}}\right)  \succ\left(  v_{i_{h+1}},v_{l_{h+1}%
}\right)
\]
and whose span is not much smaller than that of $\left(  v_{i_{1}},v_{l_{1}%
}\right)  $. Indeed, let $\left(  v_{i_{h}},v_{l_{h}}\right)  \in PR\left(
j\right)  $ be with
\begin{equation}
l_{h}-i_{h}\geq2t+5. \label{in2}%
\end{equation}

The set of all pairs $\left(  v_{i},v_{l}\right)  $ such that
\begin{equation}
\left(  v_{i},v_{l}\right)  \in PR\left(  j\right)  ,\text{ }\left(  v_{i_{h}%
},v_{l_{h}}\right)  \succ\left(  v_{i},v_{l}\right)  \label{pairs1}%
\end{equation}

is not empty - otherwise no vertex of $\left\{  v_{i_{h}+1},...,v_{l_{h}%
-1}\right\}  $ is joined to $v_{2k+1}$ and hence,
\[
l_{h}-i_{h}+1\leq\left|  M\right|  \leq t,
\]
a contradiction with (\ref{in2}). Choose a pair $\left(  v_{i_{h+1}%
},v_{l_{h+1}}\right)  $ satisfying (\ref{pairs1}) with maximal span; thus, no
vertex of
\[
\left\{  v_{i_{h}+1},...,v_{i_{h+1}-1}\right\}  \cup\left\{  v_{l_{h+1}%
+1},...,v_{l_{h}-1}\right\}
\]
is joined to $v_{2k+1}.$ Hence, we find that
\[
\left(  l_{h}-i_{h}-1\right)  -\left(  l_{h+1}-i_{h+1}+1\right)  \leq\left|
M\right|  \leq t,
\]
and thus,%
\[
l_{h+1}-i_{h+1}\geq l_{h}-i_{h}-\left(  t+2\right)  .
\]

Repeating the same argument we construct a sequence $\left(  v_{i_{h}%
},v_{l_{h}}\right)  \in PR\left(  j\right)  ,$ $h=1,...,m$ such that for every
$h=1,...,m-1,$%
\begin{align}
\left(  v_{i_{h}},v_{l_{h}}\right)   &  \succ\left(  v_{i_{h+1}},v_{l_{h+1}%
}\right)  ,\text{ }l_{h}-i_{h}\geq2t+5\nonumber\\
l_{h}-i_{h}  &  >l_{h+1}-i_{h+1}\geq l_{h}-i_{h}-\left(  t+2\right)  ,
\label{in1}%
\end{align}
and
\begin{equation}
l_{m}-i_{m}\leq2t+4. \label{in3}%
\end{equation}

Select some $h\in\left[  m-1\right]  $ and observe there are at least
$\left\lceil \left(  l_{h}-i_{h}\right)  /2\right\rceil -1$ paths of the type
$\left(  v_{2s-1},v_{2s},v_{2s+1}\right)  $ along the path $P=\left(
v_{i_{h}},...,v_{l_{h}}\right)  .$ One of these paths contain the edge
$\left(  v_{j},v_{j+1}\right)  $ and each one of the remaining can be replaced
independently by the chord joining its ends, thus shortening $P$ by $1$. In
this way we see that for every
\[
q\in\left[  \left\lfloor \frac{l_{h}-i_{h}}{2}\right\rfloor +3,l_{h}%
-i_{h}+1\right]
\]
there are $q$-reductions of $P$ that contain the edge $\left(  v_{j}%
,v_{j+1}\right)  .$ Since the ends of $P,$ and so, the ends of each of its
reductions, are joined to $v_{2k+1},$ it follows that
\[
\left[  \left\lfloor \frac{l_{h}-i_{h}}{2}\right\rfloor +4,l_{h}%
-i_{h}+2\right]  \subset L_{S}\left(  v_{j},v_{j+1}\right)  .
\]
We shall show that the shortest of these reductions of $P$ has order at most
$l_{h+1}-i_{h+1}+2.$ Indeed, assume
\[
\left\lfloor \frac{l_{h}-i_{h}}{2}\right\rfloor +3\geq l_{h+1}-i_{h+1}+2.
\]
Hence, from (\ref{in1}), we see that
\[
\left\lfloor \frac{l_{h}-i_{h}}{2}\right\rfloor +3\geq l_{h+1}-i_{h+1}+3\geq
l_{h}-i_{h}-t+1,
\]
and after simple calculations we obtain $2t+4\geq l_{h}-i_{h},$ a
contradiction with (\ref{in2}).

Therefore, for $h=1,...,m-1$ the intervals
\[
\left[  \left\lfloor \frac{l_{h}-i_{h}}{2}\right\rfloor +4,l_{h}%
-i_{h}+2\right]  ,\text{ }\left[  \left\lfloor \frac{l_{h+1}-i_{h+1}}%
{2}\right\rfloor +4,l_{h+1}-i_{h+1}+2\right]
\]
are contiguous or overlap and thus, their union is also an interval. From
(\ref{eq1}) and (\ref{in3}) we obtain%
\[
\left[  2k-d+6,2k+1\right]  \subset R_{S}\left(  v_{j},v_{j+1}\right)  ,
\]
as required.

\textbf{Case 2.}\emph{ }$x=v_{1},$\emph{ }$y=v_{2k+1}$

Observe that in the proof of the previous case we have shown that for every
\[
\left\{  v_{j},v_{j+1}\right\}  \subset\left\{  v_{1},...,v_{2k}\right\}
\]
and for every
\begin{equation}
q\in\left[  2k-d+6,2k+1\right]  , \label{incl}%
\end{equation}

there is a cycle of order $q$ of the form%
\[
\left(  v_{2k+1},v_{i},...,v_{j},v_{j+1},...,v_{i+q-2},v_{2k+1}\right)  .
\]
Applying this assertion to $\left(  v_{j},v_{j+1}\right)  =\left(  v_{1}%
,v_{2}\right)  $ we see that for every $q$ satisfying (\ref{incl}), there is a
cycle of order $q$ of the form%
\[
\left(  v_{2k+1},v_{1},v_{2},...,v_{q-1},v_{2k+1}\right)
\]
and therefore,%
\[
\left[  2k-d+6,2k+1\right]  \subset R_{S}\left(  v_{1},v_{2k+1}\right)  .
\]

\end{proof}

\begin{proof}
[\textbf{Proof of Lemma \ref{flat}}]Denote by $C$ the cycle $\left(
v_{1},...v_{2k+1},v_{1}\right)  .$ Set $t=2k-d$ and observe that from
$d_{S}\left(  v_{2k}\right)  \geq d$ and $d_{S}\left(  v_{2k+1}\right)  \geq
d,$ we have
\begin{align}
\left|  \left\{  v_{1},...,v_{2k}\right\}  \backslash\Gamma\left(
v_{2k+1}\right)  \right|   &  \leq t,\label{bnd1}\\
\left|  \left\{  v_{1},...,v_{2k-1},v_{2k+1}\right\}  \backslash\Gamma\left(
v_{2k}\right)  \right|   &  \leq t. \label{bnd2}%
\end{align}

Suppose that the distance between $x$ and $y$ along $C$ is at most $t+2$ and
let $P\left(  x,y\right)  $ be the longer $xy$-path along $C.$ Clearly,
\[
\left|  P\left(  x,y\right)  \right|  \geq2k+1-t=d+1.
\]
Hence, setting $l=\left|  P\left(  x,y\right)  \right|  $ and applying Lemma
\ref{pr1}, (i), we find that
\[
\left[  \left\lfloor \frac{l}{2}\right\rfloor +2,l\right]  \subset
R_{S}\left(  x,y\right)  ,
\]
and since for $l\geq d+2$ we have%
\[
\left\lfloor \frac{l}{2}\right\rfloor +1\leq l-\left\lceil \frac{d}%
{2}\right\rceil ,
\]
the assertion is proved in this case. So we shall hereafter assume that the
distance between $x$ and $y$ along the cycle $\left(  v_{1},...v_{2k+1}%
,v_{1}\right)  $ is at least $t+2.$

\textbf{Case 1. }$\left\{  x,y\right\}  \subset\left\{  v_{1},...,v_{2k-1}%
\right\}  $

Let $x=v_{i},$ $y=v_{j};$\ without loss of generality we assume $i<j;$ hence,
\[
j-i\geq t+2.
\]
Our first goal is to show that there exist two vertices $v_{p},v_{q}$ such
that
\begin{equation}
i<p,\text{ }p+\left(  j-i-t-2\right)  \leq q<j, \label{bnd3}%
\end{equation}

and either the edges $\left(  v_{2k},v_{p}\right)  ,$ $\left(  v_{2k+1}%
,v_{q}\right)  $ or the edges $\left(  v_{2k},v_{q}\right)  ,$ $\left(
v_{2k+1},v_{p}\right)  $ exist. Indeed, observe that the set $\left\{
v_{i+1},...,v_{j-1}\right\}  $ has at least $t+1$ members; therefore,
\[
\left\{  v_{i+1},...,v_{j-1}\right\}  \cap\Gamma\left(  v_{2k}\right)
\neq\varnothing,
\]
and%
\[
\left\{  v_{i+1},...,v_{j-1}\right\}  \cap\Gamma\left(  v_{2k+1}\right)
\neq\varnothing.
\]
Among the vertices
\[
\left\{  v_{i+1},...,v_{j-1}\right\}  \cap\left(  \Gamma\left(  v_{2k}\right)
\cup\Gamma\left(  v_{2k+1}\right)  \right)
\]
let $v_{p}$ be the one with minimal index; assume without loss of generality
that $v_{p}\in\Gamma\left(  v_{2k+1}\right)  .$ Among the vertices
\[
\left\{  v_{i+1},...,v_{j-1}\right\}  \cap\Gamma\left(  v_{2k}\right)
\]
let $v_{q}$ be the one having the maximal index. By our choice the edges
$\left(  v_{2k},v_{q}\right)  ,$ $\left(  v_{2k+1},v_{p}\right)  $ exist.
Clearly
\[
\left(  \left\{  v_{i+1},...,v_{p-1}\right\}  \cup\left\{  v_{q+1}%
,...,v_{j-1}\right\}  \right)  \cap\Gamma\left(  v_{2k}\right)  =\varnothing,
\]
implying (\ref{bnd3}).

Consider now the paths
\begin{align*}
P_{1}  &  =\left(  v_{i},v_{i-1},...,v_{1},v_{2k+1}\right)  ,\\
P_{2}  &  =\left(  v_{p},v_{p+1},...,v_{q}\right)  ,\\
P_{3}  &  =\left(  v_{2k},v_{2k-1},...,v_{j}\right)  ,
\end{align*}
and set $l_{i}=\left|  P_{i}\right|  ,$ $i=1,2,3.$ The concatenation
$Q=\left(  P_{1},P_{2},P_{3}\right)  $ is a $v_{i}v_{j}$-path with
\[
\left|  Q\right|  =l_{1}+l_{2}+l_{3}+2=2k+1-\left(  j-i\right)  +\left(
q-p\right)  +2.
\]
Set $l=\left|  Q\right|  ;$ from (\ref{bnd3}) we obtain%
\[
l\geq2k+1-\left(  j-i\right)  +\left(  q-p\right)  +2\geq2k+1-t=d+1.
\]
Applying Lemma \ref{pr1}, part (i), to each one of the paths $P_{1}%
,P_{2},P_{3},$ we see that $Q$ has a $q$-reduction for every
\[
q\in\left[  l-\left\lfloor \frac{l_{1}}{2}\right\rfloor -\left\lfloor
\frac{l_{2}}{2}\right\rfloor -\left\lfloor \frac{l_{3}}{2}\right\rfloor
+3,l\right]  .
\]
In view of
\[
\left\lfloor \frac{l_{1}}{2}\right\rfloor +\left\lfloor \frac{l_{2}}%
{2}\right\rfloor +\left\lfloor \frac{l_{3}}{2}\right\rfloor \geq\left\lfloor
\frac{l_{1}+l_{2}+l_{3}}{2}\right\rfloor -1=\left\lfloor \frac{l}%
{2}\right\rfloor -2\geq\left\lceil \frac{d}{2}\right\rceil -2,
\]
the assertion follows.

\textbf{Case }$2$\textbf{. }$x=v_{2k+1},$\emph{ }$y\in\left\{  v_{1}%
,...,v_{2k-1}\right\}  $

Let $y=v_{j};$ since the $xy$-distance along $C$ is at least $t+2,$ we have
$j\geq t+2,$ and thus%
\[
M=\left\{  v_{1},...,v_{j-1}\right\}  \cap\Gamma\left(  v_{2k}\right)
\neq\varnothing,
\]
Among the vertices of $M$ let $v_{q}$ be the one having the maximal index.
Clearly,
\[
\left\{  v_{q+1},...,v_{j-1}\right\}  \cap\Gamma\left(  v_{2k}\right)
=\varnothing,
\]
and thus,
\begin{equation}
j-q-1\leq t \label{bnd4}%
\end{equation}

Consider now the paths
\begin{align*}
P_{1}  &  =\left(  v_{2k+1},v_{1},...,v_{q}\right)  ,\\
P_{2}  &  =\left(  v_{2k},v_{2k-1},...,v_{j}\right)
\end{align*}
and set $l_{i}=\left|  P_{i}\right|  ,$ $i=1,2.$ The concatenation $Q=\left(
P_{1},P_{2}\right)  $ is a $v_{2k+1}v_{j}$-path with
\[
\left|  Q\right|  =l_{1}+l_{2}+1=2k+1-j+q+1.
\]
Set $l=\left|  Q\right|  ;$ from (\ref{bnd4}) we obtain%
\[
l\geq2k+1-j+q+1\geq2k+1-t=d+1.
\]
Applying Lemma \ref{pr1}, part (i), to each one of the paths $P_{1},P_{2},$ we
see that $Q$ has a $q$-reduction for every
\[
q\in\left[  l-\left\lfloor \frac{l_{1}}{2}\right\rfloor -\left\lfloor
\frac{l_{2}}{2}\right\rfloor +2,l\right]  .
\]
In view of
\[
\left\lfloor \frac{l_{1}}{2}\right\rfloor +\left\lfloor \frac{l_{2}}%
{2}\right\rfloor \geq\left\lfloor \frac{l_{1}+l_{2}}{2}\right\rfloor
-1=\left\lfloor \frac{l-1}{2}\right\rfloor -1\geq\left\lceil \frac{d}%
{2}\right\rceil -2,
\]
the assertion follows.

\textbf{Case }$3.$\emph{ }$y=v_{2k},$\emph{ }$x\in\left\{  v_{1}%
,...,v_{2k-1}\right\}  $

This case is symmetric to the previous. Setting $x=v_{j},$ we find a vertex
\[
v_{q}\in\left\{  v_{j+1},...,v_{2k}\right\}  \cap\Gamma\left(  v_{2k+1}%
\right)  ,
\]
then consider the paths%
\begin{align*}
P_{1}  &  =\left(  v_{2k},v_{2k-1},...,v_{q}\right)  ,\\
P_{2}  &  =\left(  v_{2k+1},v_{1},v_{2},...,v_{j}\right)  ,
\end{align*}

and find $xy$-paths of proper order among the reductions of the concatenation
$\left(  P_{1},P_{2}\right)  .$
\end{proof}

\begin{proof}
[\textbf{Proof of Lemma \ref{le1.3}}]Applying Lemma \ref{pr1}, part (i), to
the path
\[
P=\left(  v_{1},...,v_{2k+1}\right)
\]
we see that $S$ contains cycles of order $q$ for every $s\in\left[
k+1,2k+1\right]  $ and the proof is completed under the assumption $k+1\leq4r$.

Assume now that $k\geq4r;$ let
\[
S_{1}=\left\{  v_{1},...,v_{4r-1}\right\}  ,\text{ }S_{2}=\left\{
v_{4r},...,v_{2k+1}\right\}  .
\]
We have shown that for every $q\in\left[  2r,4r-1\right]  $ there is
a\ $q$-reduction of the path $\left(  v_{1},...,v_{4r-1}\right)  .$ On the
other hand the order of the path
\[
P=\left(  v_{4r},...,v_{2k+1}\right)
\]
is $\left(  2k-4r+2\right)  ,$ and applying the Chopping Lemma, we see that
for every interval
\[
I\subset\left[  2k-2r+1\right]  ,\text{ }\left|  I\right|  =2r,
\]
there is a $q$-reduction of $P$ for some $q\in I$. Applying the Collating
Lemma to the graph $S$ with the partition
\[
V\left(  S\right)  =S_{1}\cup S_{2}%
\]
and the edges $\left(  v_{1},v_{2k+1}\right)  $ and $\left(  v_{4r+1}%
,v_{4r+2}\right)  $ it follows that $S$ contains cycles of order $q$ for
every
\[
q\in\left[  4r,2k+1\right]  .
\]

\end{proof}

\subsection{Concluding remarks and open problems}

There is a much simpler proof of (\ref{maineq}) under the assumption
$p\geq8r+7.$ Actually, except for Lemma \ref{le1.3}, our methods are good
enough to prove (\ref{maineq}) for $p\geq3r+9,$ and\ it seems that with some
additional refinement it is possible to prove (\ref{maineq}) for
\[
p\geq2r+o\left(  r\right)  .
\]
The following conjecture, however, looks more challenging.

\begin{conjecture}
For every $k$ there exists $r_{0}=r_{0}\left(  k\right)  $ such that for
$r>r_{0}$ and $p>r^{1/k},$%
\[
r\left(  C_{p},K_{r}\right)  =\left(  p-1\right)  \left(  r-1\right)  +1.
\]

\end{conjecture}

There are known Ramsey numbers $r\left(  C_{p},K_{r}\right)  $ for $p<r$ -
Jayawardene\ and Rousseau found that $r\left(  C_{4},K_{6}\right)  =18$ in,
\cite{JaRo97} and $r\left(  C_{5},K_{6}\right)  =21$ in \cite{JaRo00};
Schiermeyer found that $r\left(  C_{5},K_{7}\right)  =25$ in \cite{Schr03}.
These values, although very few, give some hope that the conjecture might be true.

\textbf{Acknowledgement }The author is grateful to Cecil Rousseau and to Dick
Schelp for the many delightful hours spent discussing Ramsey problems and the
above problem in particular. The help and advice of Cecil Rousseau were in
every respect invaluable. Finally, B\'{e}la Bollob\'{a}s suggested many
corrections and improvements of the manuscript.


\begin{thebibliography}{99}                                                                                               %


\bibitem {BJM}B. Bollob\'{a}s, C. Jayawardene, J. Yang, Y. Huang, C. C.
Rousseau and K. Zhang, On a conjecture involving cycle-complete graph Ramsey
numbers, \emph{Australas. J. Combin. }\textbf{22} (2000), 63--71.

\bibitem {Bol}B. Bollob\'{a}s, \emph{Modern graph theory,} Graduate Texts in
Mathematics, 184, Springer-Verlag, New York (1998), xiv+394 pp.

\bibitem {Bon}J. A. Bondy, Basic graph theory: paths and circuits,
\emph{Handbook of combinatorics}, Vol. \textbf{1}, Elsevier, Amsterdam, 1995,
pp. 3--110.

\bibitem {BEFRS}S. Burr, P. Erd\H{o}s, R. J. Faudree, C. C. Rousseau and R. H.
Schelp, An extremal problem in generalized Ramsey theory, \emph{Ars
Combinatoria,} \textbf{10} (1980), 193--203.

\bibitem {BoEr}J. A. Bondy and P. Erd\"{o}s, Ramsey numbers for cycles in
graphs. \emph{J. Comb. Theory Ser. B} \textbf{14} (1973), 46--54.

\bibitem {Dir}G. A. Dirac, Some theorems on abstract graphs,
$\emph{Proc.London}$ $\emph{Math.}$ $\emph{Soc.}$\textbf{2} (1952), 69--81.

\bibitem {ErGa}P. Erd\H{o}s and T. Gallai, On maximal paths and circuits of
graphs, \emph{Acta Math. Acad. Sci. Hungar.} \textbf{10} 1959, 337--356.

\bibitem {EFRS}P. Erd\H{o}s, R. J. Faudree, C. C. Rousseau and R. H. Schelp,
On cycle-complete graph Ramsey numbers, \emph{J. Graph Theory }\textbf{2}
(1978), 53--64.

\bibitem {JaRo97}Ch. Jayawardene and C.C. Rousseau, The Ramsey number for a
quadrilateral vs. a complete graph on six vertices, \emph{Congr. Numer.}
\textbf{123} (1997), 97--108.

\bibitem {JaRo00}Ch. Jayawardene and C.C. Rousseau, The Ramsey number for a
cycle of length five vs. a complete graph of order six. \emph{J. Graph Theory}
\textbf{35} (2000), 99--108.

\bibitem {Rou}C. C. Rousseau, Asymptotic bounds for Ramsey Numbers, preprint.

\bibitem {Schr}I. Schiermeyer, All Cycle-Complete graph Ramsey Numbers
$R(C_{m},K_{6})$, \emph{J. Graph Theory} \textbf{44} (2003), 251-260.

\bibitem {Schr03}I. Schiermeyer, The Cycle-Complete graph Ramsey Numbers
$R(C_{5},K_{7})$, preprint, 2003.

\bibitem {YHZ}J. Yang, Y. Huang and K. Zhang, The value of the Ramsey number
$R\left(  C_{n},K_{4}\right)  $ is $3\left(  n-1\right)  +1$ $\left(
n\geq4\right)  $, \emph{Australas. J. Combin.} \textbf{20} (1999), 205--206.

E-mail address: \textit{vnikifrv@memphis.edu}
\end{thebibliography}
\end{document}